\documentclass{article}

\usepackage{microtype}
\usepackage{graphicx}
\usepackage{subfigure}
\usepackage{booktabs} 

\usepackage{hyperref}

\usepackage[accepted]{icml2024}

\usepackage{amsmath}
\usepackage{amssymb}
\usepackage{mathtools}
\usepackage{amsthm}
\usepackage{multirow}
\usepackage{tikz}
\usetikzlibrary{positioning}
\usetikzlibrary{arrows.meta}
\usepackage[capitalize,noabbrev]{cleveref}

\theoremstyle{plain}
\newtheorem{theorem}{Theorem}[section]

\theoremstyle{definition}

\theoremstyle{remark}

\usepackage[textsize=tiny]{todonotes}

\icmltitlerunning{Neural operators meet conjugate gradients: The FCG-NO method for efficient PDE solving}

\begin{document}

\twocolumn[
\icmltitle{Neural operators meet conjugate gradients: The FCG-NO method for efficient PDE solving}



\icmlsetsymbol{equal}{*}

\begin{icmlauthorlist}
\icmlauthor{Alexander Rudikov}{Sk}
\icmlauthor{Vladimir Fanaskov}{Sk}
\icmlauthor{Ekaterina Muravleva}{Sk}
\icmlauthor{Yuri M. Laevsky}{ICMMG}
\icmlauthor{Ivan Oseledets}{AIRI,Sk}
\end{icmlauthorlist}

\icmlaffiliation{Sk}{Skolkovo Institute of Science and Technology}
\icmlaffiliation{AIRI}{Artificial Intelligence
Research Institute}
\icmlaffiliation{ICMMG}{Institute of Computational Mathematics and Mathematical Geophysics SB RAS}

\icmlcorrespondingauthor{Alexander Rudikov}{A.Rudikov@skoltech.ru}

\icmlkeywords{Machine Learning,ICML}

\vskip 0.3in
]



\printAffiliationsAndNotice{}  

\begin{abstract}
Deep learning solvers for partial differential equations typically have limited accuracy. We propose to overcome this problem by using them as preconditioners. More specifically, we apply discretization-invariant neural operators to learn preconditioners for the flexible conjugate gradient method (FCG). Architecture paired with novel loss function and training scheme allows for learning efficient preconditioners that can be used across different resolutions. On the theoretical side, FCG theory allows us to safely use nonlinear preconditioners that can be applied in $O(N)$ operations without constraining the form of the preconditioners matrix. To justify learning scheme components (the loss function and the way training data is collected) we perform several ablation studies. Numerical results indicate that our approach favorably compares with classical preconditioners and allows to reuse of preconditioners learned for lower resolution to the higher resolution data.
\end{abstract}

\section{Introduction}
The recent surge of interest in learning solution operators and surrogates for partial differential equations (PDEs) leads to several new approaches and architectures \cite{azzizadenesheli2023neural}, \cite{karniadakis2021physics}, \cite{kovachki2021neural}. Most notable, combination of functional methods with deep learning (spectral convolutions \cite{rippel2015spectral}, FNO \cite{li2020fourier}); operator-valued kernels \cite{kadri2016operator}, \cite{batlle2024kernel}; random feature model in Banach space; \cite{nelsen2021random} architectures based on universal approximation for operators (DeepONet \cite{lu2019deeponet} and its variants).

According to the large-scale benchmarks, typical accuracy of such methods is about $0.1\%-1\%$ relative $L_2$ norm \cite{takamoto2022pdebench}, \cite{lu2022comprehensive}, \cite{de2022cost}. On the other hand, classical numerical methods for PDEs are usually consistent, so they can reach arbitrary accuracy on a sufficiently fine grid.

\begin{figure}[h]
  \centering
    \includegraphics[scale=0.67]{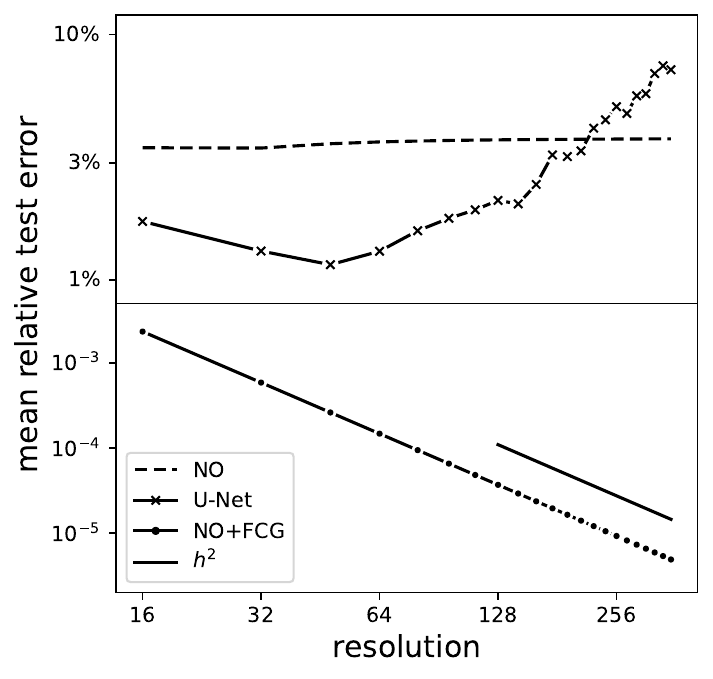}
  \caption{Comparison of accuracy for three approaches: U-Net -- classical deep learning architecture, NO -- neural operator, NO+FCG -- hybrid approach advocated in the present article. Due to the prominent difference between accuracies, two distinct scales are used in the $y$ axis. One can observe that owing to the finite receptive field the performance of U-Net deteriorates with the increase of resolution. The neural operator provides the same accuracy with resolution increase -- this is a highly-praised ``discretization invariance.'' For NO+FCG, the error decreases with the increase of resolution in the same way as for classical numerical methods.}
  \label{fig:superresolution}
\end{figure}

Is it possible to combine the consistency with the efficiency of deep learning methods? One way around this is to build a hybrid system, that is, to replace some parts of the classical method with deep learning components (e.g., \cite{bar2019learning}, \cite{kochkov2021machine}, \cite{greenfeld2019learning}). In the present contribution we focus on elliptic boundary value problems and show how to utilize neural operator \cite{li2020fourier} (NO) with flexible conjugate gradient method \cite{notay2000flexible} (FCG). As one can see on \cref{fig:superresolution}, our approach indeed allows one to retain consistency which leads to much better error on grids with higher resolutions.

We are not the first ones, who propose to learn preconditioners. The previous contributions include at least \cite{hsieh2019learning}, \cite{zhang2022hybrid}, \cite{cui2022fourier}, \cite{hausner2023neural}, \cite{li2023learning}. However, our approach is the first one with three unique properties: (i) -- we learn nonlinear operator with asymptotic complexity $O(N)$, which allows us not to prescribe a structure of the preconditioners matrix to have an efficient matrix-vector product, (ii) -- our nonlinear operator is discretization-invariant so it can be applied at different resolutions, (iii) -- we have convergence guarantees based on FCG theory built in \cite{notay2000flexible}. We provide more details on comparison with other approaches in \cref{section:related research}.

To summarize, our contributions are:
\begin{enumerate}
    \item We propose to train a neural operator as a nonlinear preconditioner for the flexible conjugate gradient method. We demonstrate that neural operators can be trained on lower resolutions and serve as an efficient preconditioner for higher resolutions.
    \item We provide a novel learning scheme that involves drawing random vectors from the Krylov subspace. Ablation study shows that abandoning Krylov subspace and using classical random right-hand sides seriously impair training.
    \item We put forward a novel loss based on energy norm that comes with convergence guarantees and significantly outperforms $L_2$ loss usually used for training. 
\end{enumerate}

\label{section:introduction}

\section{Neural operator as a preconditioner}
Here, we introduce the PDE that we solve and provide a brief overview of discretization, iterative schemes, and linear and nonlinear preconditioning. After that, we outline the construction of the loss function, neural operator, and training strategy.

\subsection{Elliptic equations}
A family of PDEs we consider in this article has a form
\begin{equation}
    \label{eq:general_elliptic_equation}
    \begin{split}
    -&\sum_{ij=1}^{2}\frac{\partial}{\partial x_i}\left(a(x) \frac{\partial u(x)}{\partial x_j} \right) = f(x)\\
    &x \in \Gamma \equiv (0,~1)^2, ~u(
    x)\big|_{x \in \partial{\Gamma}} = 0,
    \end{split}
\end{equation}
where $\partial{\Gamma}$ is a boundary of the unit hypercube $\Gamma$, and $a(x)\geq\epsilon>0$. Applying the finite element method or finite difference method (FDM) to \eqref{eq:general_elliptic_equation} we can obtain a large sparse linear algebraic system with a symmetric positive definite matrix:
\begin{equation}
    \label{eq:linear_system}
    Au = f;\,u,\,f \in \mathbb{R}^{n},\,A\in\mathbb{R}^{n\times n}, A \succ 0.
\end{equation}
The system can be solved in $O\left(n^3\right)$ operations by Gauss elimination, but for large $n^{3}$ one often resorts to iterative techniques that can better leverage the special sparsity stricture of matrix $A$.

\subsection{Preconditioning of linear system}
One of the commonly used iterative methods for solving \eqref{eq:linear_system} is conjugate gradient (CG). However, in practice, the CG often converges slowly, and the convergence also depends on the distribution of the eigenvalues and the initial residual. One way to enhance convergence is to precondition the system to improve its spectral properties. The modified linear system becomes
\begin{equation}
    \label{eq:def_of_preconditioner}
    B^{-1} A u = B^{-1} f,
\end{equation}
where $B \in \mathbb{R}^{n \times n}$ is 
called the preconditioner.\footnote{Note that formally for $A\succ 0$ one should use $B^{-\frac{1}{2}} A B^{-\frac{1}{2}}$ to preserve properties of the system. But since $B^{-\frac{1}{2}} A B^{-\frac{1}{2}}$ and $B^{-1}A$ are spectrally equivalent one can rewrite conjugate gradient algorithm to work directly with $B^{-1}A$ using $B$ inner product. See Section 9.2 from \cite{saad2003iterative} for more details.}
To be successful, preconditioner $B$ should have several properties: (i) it should improve spectral properties of $A$, e.g., the condition number of $BA$ is much smaller than the condition number of $A$; (ii) $B$ easily invertible, i.e., $Bg=r$ is cheap to solve;  (iii) for preconditioner conjugate gradient (PCG) \cite{axelsson1996iterative} one need to ensure $B$ is symmetric positive definite. All these requirements greatly complicate the construction of preconditioners. Fortunately, with a slight decrease in numerical efficiency, one can switch to a less restrictive set of nonlinear preconditioners.


\begin{algorithm}[tb]
   \caption{Flexible Conjugate Gradients}
   \label{alg:algorithm_1}
    \begin{algorithmic}
        \STATE {\bfseries Input:} $A$, $\mathcal{B}$, $f$, $m_{\text{max}} > 0$, $\text{iter}$.
        \STATE {\bfseries Ensure:} $u_{\text{iter}}$,\,$r_{\text{iter}}$.
        \STATE Initialize $u_0 \leftarrow \mathcal{N}(0,\,1) \in \mathbb{R}^n$, $r_0 \leftarrow f - Au_0 \in \mathbb{R}^n$.
        \FOR{$i=0$ {\bfseries to} $\text{iter} - 1$}
            \STATE $w_i \leftarrow \mathcal{B}(r_i)$
            \STATE $m_i \leftarrow \min(i,\,\max(1,\mod(i,\,m_{\text{max}} + 1)))$
            \STATE $p_i \leftarrow w_i - \sum_{k=i-m_i}^{i-1}\dfrac{\left(w_i, s_k\right)}{\left(p_k, s_k\right)}~p_k$
            \STATE $s_i \leftarrow A p_i$
            \STATE $u_{i+1} \leftarrow u_i + \dfrac{\left(p_i, r_i\right)}{\left(p_i, s_i\right)}~p_i$
            \STATE $r_{i+1} \leftarrow r_i - \dfrac{\left(p_i, r_i\right)}{\left(p_i, s_i\right)}~s_i$
        \ENDFOR
    \end{algorithmic}
\end{algorithm}

\subsection{Nonlinear preconditioners}
A projection iterative method suitable for more general preconditioners appeared in \cite{notay2000flexible}. The algorithm is called Flexible Conjugate Gradient (FCG) and is given in \cref{alg:algorithm_1}. It uses many vectors instead of one so it is less effective than CG from a computational perspective. But in return, one can consider $B$ in \cref{eq:def_of_preconditioner} to be the nonlinear operator. The only technical requirement is given in the following result from \cite{notay2000flexible}:
\begin{theorem}
    \label{thm:theorem_1}
    Let $A,\,B \in \mathbb{R}^{n \times n}$ be  symmetric positive definite matrices and $\mathcal{B}: \mathbb{R}^{n} \rightarrow \mathbb{R}^{n}$. Let $f,~u_0$ be the vectors of $\mathbb{R}^n$, and let $\big\{r_i\big\}_{i=0,1,\ldots},\,\big\{p_i\big\}_{i=0,1,\ldots},\,\big\{u_i\big\}_{i=1,2,\ldots}$ be the sequences of vectors generated by applying \cref{alg:algorithm_1} to $A,~\mathcal{B},~f$, and $u_0$ with some given sequences of non-negative integer parameters $\big\{m_i\big\}_{i=0,1,\ldots}$.

    If, for any $i$,
    \begin{equation}
        \label{eq:theorem}
        \dfrac{\big \Vert \mathcal{B}(r_i) - B^{-1}r_i \big \Vert_B}{\big \Vert B^{-1}r_i \big \Vert_B} \leqslant \varepsilon_i < 1, 
    \end{equation}
    then
    \begin{equation*}
        \dfrac{\big\Vert u - u_{i+1}\big\Vert_A}{\big\Vert u - u_i\big\Vert_A} \leqslant \dfrac{\kappa\big(B^{-1}A\big)\cdot \gamma_i - 1}{\kappa\big(B^{-1}A\big)\cdot \gamma_i + 1},
    \end{equation*}
    where $\gamma_i = \dfrac{1 + \varepsilon_i}{1 - \varepsilon_i} \cdot \dfrac{\left(1 + \varepsilon^2_i\right)^2}{\left(1 - \varepsilon_i^2\right)}$, and $~\big\Vert u \big\Vert_A = \sqrt{\left(u, A u\right)}$.
\end{theorem}
In plain English, \cref{thm:theorem_1} states that as long as nonlinear preconditioner $\mathcal{B}$ is close enough to the linear preconditioner $B$ (in a sense of equation~\eqref{eq:theorem}), FCG converges with roughly the same speed as CG with $B$ taken as a preconditioner. This powerful result allows us to use a nonlinear neural network as a preconditioner and simultaneously provide a loss, suitable for training this neural network. Both these points are explained in more detail in the next section.

\subsection{Learning scheme}
\subsubsection{Loss functions}
The principal idea of the article is to select a family of neural networks with weights $\theta$ as nonlinear preconditioner $\mathcal{B}(r; \theta)$ for FCG described in   \cref{alg:algorithm_1}. \cref{thm:theorem_1} can be interpreted in terms of optimization target. Namely, we can find parameters from the loss
\begin{equation}
    \label{eq:unrealistic_loss}
    L_{\text{opt}}(\theta) = \max_{r}\dfrac{\big \Vert \mathcal{B}(r; \theta) - A^{-1}r \big \Vert_A}{\big \Vert A^{-1}r \big \Vert_A},\text{ s.t. } \left\|r\right\|_2 = 1.
\end{equation}
If we were able to minimize \cref{eq:unrealistic_loss} with respect to $\theta$ we would find a neural network that provides a nonlinear preconditioner as close to $A^{-1}$ as possible. Unfortunately \eqref{eq:unrealistic_loss} is a hard minimax problem so we relax it to
\begin{equation}
    \label{eq:relaxed_loss}
    L_{\text{mean}}(\theta) = \mathbb{E}_{r}\dfrac{\big \Vert \mathcal{B}(r; \theta) - A^{-1}r \big \Vert_A}{\big \Vert A^{-1}r \big \Vert_A}.
\end{equation}
To tie this loss function to the usual operator learning framework described in \cite{li2020fourier} we suppose that parameters of PDE \eqref{eq:general_elliptic_equation} are drawn from some distributions in Banach space $a\sim p_a$ and $f\sim p_f$ after that the final loss, which we call Notay loss after the author of \cite{notay2000flexible}, becomes
\begin{equation}
    \label{eq:final_loss}
    L_{\text{Notay}}(\theta) = \mathbb{E}_{r, a, f}\dfrac{\big \Vert \mathcal{B}(r; \theta) - A^{-1}r \big \Vert_A}{\big \Vert A^{-1}r \big \Vert_A},
\end{equation}
where $A$ depends on $a$ via discretization, $r$ depends on $f$ and the distribution of initial guess $u_0\sim \mathcal{N}(0, I)$. Even though $r$ is not independent from other variables we still have freedom in the choice of distribution for $r$.

The simplest case is to compute the residual from $u_0$ and $f$. However, one expects that distribution for $r$ will drift from $r = f - Au_0,\,f\sim p_f,\,u_0\sim \mathcal{N}(0, I)$ since both CG and FCG seek a solution in the Krylov subspace
\begin{equation}
    \mathcal{K}_{m}(A, r_0) = \text{Span}\Big\{r_0, Ar_0, \ldots, A^{m-1}r_0\Big\}.
\end{equation}
So it is natural to consider a more general distribution for $r$ drawing vectors from Krylov subspaces using distributions for $u_0$ and $f$
\begin{equation}
    \label{eq:rk_distribution}
    r \sim p_{\mathcal{K}_{m}}(r) \Leftrightarrow r \in \mathcal{K}_{m}(A, r_0),\,f\sim p_f,\,u_0\sim \mathcal{N}(0, I).
\end{equation}
In other words, to draw residuals from $p_{\mathcal{K}_{m}}(r)$ we draw $u_0$ and $f$ and run CG for $m$ iterations.

Of course a valid alternative to \cref{eq:final_loss} is an ordinary $L_2$ loss function
\begin{equation}
    \label{eq:L2_loss}
    L_{2}(\theta) = \mathbb{E}_{r, a, f}\dfrac{\big \Vert \mathcal{B}(r; \theta) - A^{-1}r \big \Vert_2}{\big \Vert A^{-1}r \big \Vert_2}.
\end{equation}
Since $L_2$ loss \eqref{eq:L2_loss} is commonly used to learn the solution operator for the elliptic equation, we will also try to use it to learn the preconditioner.

\subsubsection{Neural operator}
In principle, losses \eqref{eq:final_loss} and \eqref{eq:L2_loss} are suitable for arbitrary neural networks. However, it is desirable to learn a preconditioner for one low-resolution discretization and apply it to discretizations with higher resolutions. Models uniquely suited for that task are neural operators that are known to be consistent under the change of resolution. Following \cite{li2020fourier}, \cite{fanaskov2022spectral} we consider spectral neural operators (SNO) with linear integral kernels $$u\leftarrow \int dx A_{ij} p_j(x) \left(p_i, u\right),$$ where $p_j(x)$ are orthogonal or trigonometric polynomials.

Layers with linear integral kernels can be described in terms of analysis and synthesis matrices. To describe these matrices consider a finite set of orthogonal polynomials (with weight function $w$) $p_n(x), ~n = 0, \ldots, N-1$. The synthesis operator is 
\[
f(x) = \sum_{i=0}^{N-1}p_i(x)c_i ~\Longleftrightarrow ~f = \mathcal{S}c,
\]
where $c$ is a vector with expansion coefficients $c_i, ~i=0, \ldots, N-1$ defined by 
\[
    c_i = \dfrac{1}{h_i}\int p_i(x) f(x) w dx,\, h_i = \int p_i^2(x) w dx.
\]
When the Gauss quadrature \cite{golub1969calculation, gautschi2004orthogonal, trefethen2008gauss} is used for the integral evaluation, it is easy to show that the analysis operator has a form
$$c = \mathcal{A}f = D^{-1} \mathcal{S}^\top \big(f(x) \odot w\big),$$
where $$D = \text{diag} \Big(p_{0}^\top\big(p_0 \odot w\big), \ldots, p_{N-1}^\top\big(p_{N-1} \odot w\big)\Big).$$ So, the whole integral kernel has a form $u\leftarrow \mathcal{S} L \mathcal{A}u$ where $L$ is a linear operator on $\mathbb{R}^{N}$. For example, Fourier neural operator \cite{li2020fourier} employees FFT as analysis and inverse FFT as synthesis; neural operators from \cite{fanaskov2022spectral} utilizes DCT-II to work with Chebyshev polynomials; wavelet neural operator \cite{tripura2022wavelet} applies FWT and inverse FWT.

\label{section:preconditioner}

\label{section:scheme}

\section{Experiments}
Here, we provide proof of the effectiveness of the proposed approach. For that purpose, we design several experiments in $D=2$. We start with a description of the training details.

\subsection{Training details}

\begin{figure*}[!h]
\normalsize
\centering
    \begin{center}
        \moveright  10pt \hbox{\input{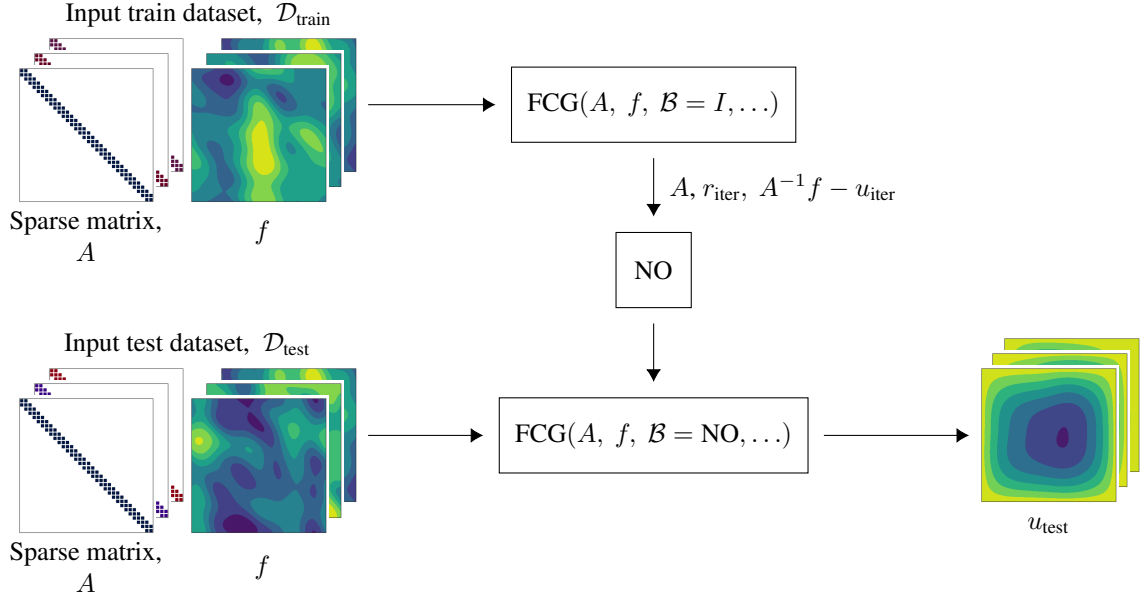}}
    \end{center}
\caption{The full scheme of the proposed approach: starts from the input train dataset, $\mathcal{D}_{\text{train}} = \big(A, f\big)$. \textbf{(a)} Submit $\mathcal{D}_{\text{train}}$ to the CG (FCG with $\mathcal{B} = I$). \textbf{(b)} Train the NO on the FCG output: $A, u_{\text{iter}}, r_{\text{iter}}$. \textbf{(c)} Apply the FCG with $\mathcal{B} = \text{NO}$ with the test dataset, $\mathcal{D}_{\text{test}}$. \textbf{(d)} Output $u_{\text{test}}$.}
  \label{fig:algorithm}
\vspace*{4pt}
\end{figure*}

According to \eqref{eq:rk_distribution}, the residuals should be derived from $m$ iterations of CG. Thus, we propose the following training scheme \cref{fig:algorithm}. Following this scheme, we start with the generation of train $\mathcal{D}_{\text{train}}$ and test $\mathcal{D}_{\text{test}}$ datasets that consist of sparse matrices $A$ and the right-hand sides $f$ of \eqref{eq:linear_system}. Then, we use $m$ iterations of the CG (not necessarily until complete convergence) on the training dataset. The CG output: $u_{\text{iter}},\,r_{\text{iter}}$ which is in the Krylov subspace, is used for the training of the neural operator (NO). The trained NO is used as the preconditioner for the FCG algorithm applied on the test dataset $\mathcal{D}_{\text{test}}$.

To accelerate the training of NO, we eliminate the need to invert the matrix $A$ in the loss function \eqref{eq:final_loss}. To accomplish this, we use the error at iteration $i$
\begin{equation*}
    e_i = u_{\text{exact}} - u_i = A^{-1}r_i,
\end{equation*}
where $u_{\text{exact}} = A^{-1} f$. After that, we can rewrite Notay loss \eqref{eq:final_loss} in the form
\begin{equation}
    \label{eq:notay_loss_with_error}
    L_{\text{Notay}}(\theta) = \mathbb{E}_{r, a, f} \dfrac{\big \Vert \mathcal{B}(r; \theta) - e \big \Vert_A}{\big \Vert e \big \Vert_A},
\end{equation}
and $L_2$-loss in the form
\begin{equation}
    \label{eq:l2_loss_with_error}
    L(\theta) = \mathbb{E}_{r, a, f} \dfrac{\big \Vert \mathcal{B}(r; \theta) - e \big \Vert_2}{\big \Vert e \big \Vert_2}.
\end{equation}

\subsection{Datasets}
To verify the proposed approach, we use two versions of the elliptic equation \eqref{eq:general_elliptic_equation} in the domain $\Gamma \equiv (0, 1)^{D}$.

First, we define random trigonometric polynomials
\begin{equation*}
    \label{eq:random_trig_2D}
    \mathcal{P}\big(N_1, N_2, \alpha\big) = \Big\{f(x) = \mathcal{R}\left(g(x)\right):\mathcal{R}(c) \simeq \mathcal{N}(0, I)\Big\},
\end{equation*}
where 
\begin{equation*}
    g(x) =\sum_{m=0}^{N_1}\sum_{n=0}^{N_2}\frac{c_{mn}\exp\left(2\pi i(mx_1 + nx_2)\right)}{(1+m+n)^\alpha}.
\end{equation*}

For the first dataset, called Poisson, we use random trigonometric
polynomials for $f$:
\begin{equation}
    \label{eq:poisson}
    a(x) = I;\,f(x) \simeq \mathcal{P}(5, 5, 2).
\end{equation}

For the second one, called Diffusion, we use trigonometric polynomials for both $a$ and $f$:
\begin{equation}
    \label{eq:diffusion}
    a(x)  \simeq \mathcal{P}(5, 5, 2) + 10;\,f(x) \simeq \mathcal{P}(5, 5, 2).
\end{equation}

\subsection{Ablation study}
In this article, we have proposed three new components:
\begin{enumerate}
    \item The nonlinear preconditioner $\mathcal{B}$ in the form of a NO in the \cref{alg:algorithm_1}; 
    \item The loss function derived from the \cref{thm:theorem_1} for training this NO; 
    \item The novel sampling strategy of a train dataset from the Krylov subspace.
\end{enumerate}

Next, we present the results of the ablation study for all these components.

\subsubsection{FCG vs. classical techniques}
\begin{table*}[!ht]
    \centering
    \begin{tabular}{@{}lcccccccc@{}}
        \toprule
        Dataset & grid & NO+FCG & CG & Jacobi(4) & GS(1) & GS(4) & ILU(1) & ILU(8) \\\midrule
        \multirow{3}{*}{Poisson} & $32$ & $\boldsymbol{9}$ & $74$ & $30$ & $27$ & $13$ & $69$ & $27$\\
        & $64$ & $\boldsymbol{14}$ & $130$ & $58$ & $47$ & $23$ & $110$ & $77$\\
        & $128$ & $\boldsymbol{20}$ & $216$ & $112$ & $78$ & $37$ & $185$ & $128$\\\cmidrule(r){2-9} 
        \multirow{3}{*}{Diffusion} & $32$ & $\boldsymbol{9}$ & $75$ & $32$ & $27$ & $13$ & $69$ & $27$\\
        & $64$ & $\boldsymbol{14}$ & $132$ & $61$ & $46$ & $22$ & $110$ & $74$\\
        & $128$ &  $\boldsymbol{19}$ & $215$ & $115$ & $78$ & $36$ & $177$ & $128$\\
        \bottomrule 
    \end{tabular}
    \caption{Comparison of FCG with classical preconditioning techniques. The table contains the first iteration number $i$ such that $\left\|r_{i}\right\|_2\big/\left\|r_{0}\right\|_2 \leq 10^{-6}$ for different resolutions.}
    \label{table:FCG_vs_classics}
\end{table*}

We perform several experiments to compare learned preconditioners with classical approaches. More specifically, we use Jacobi, symmetric Gauss-Seidel \cite{saad2003iterative} and incomplete LU (ILU) preconditioner implemented in SuperLU library \cite{demmel1999superlu}.

Recall, that Jacobi iteration reads
\begin{equation}
    \label{eq:Jacobi}
    x^{n+1} = x^{n} + D(A)^{-1}(b - Ax^{n}),
\end{equation}
where $D(A)$ is a diagonal part of matrix $A$. Symmetric Gauss-Seidel iteration reads
\begin{equation}
    \begin{split}
        &x^{n+1/2} = x^{n} + L(A)^{-1}(b - Ax^{n}),\\
        &x^{n+1} = x^{n+1/2} + U(A)^{-1}(b - Ax^{n+1/2}),
    \end{split}
\end{equation}
where $L(A)$ and $U(A)$ are lower and upper triangular parts of $A$.

When classical iteration are used to approximate $A^{-1}z$ they are used for $k$ steps (e.g., Jacobi$(4)$ means $4$ iteration of the form \eqref{eq:Jacobi}) with $x^{0} = 0$ and $b = z$.

ILU is a LU decomposition of sparse matrix $A$ with control of fill-in. For the details on the particular version used, we refer to \cite{demmel1999superlu}. In \cref{table:FCG_vs_classics} ILU$(k)$ number $k$ corresponds to parameter \textit{filling\_factor} in documentation. It prescribes the number of nonzero elements in the preconditioner. For example, for $5$-point stencil discretization of the Poisson equation, the preconditioner formed with ILU$(1)$ and ILU$(8)$ have $\simeq1.5\times$ and $\simeq6.7\times$ more nonzero elements that original matrix.
 
The results appear in \cref{table:FCG_vs_classics}. As one can see, the learned preconditioner allows for faster convergence than standard approaches on all resolutions tested. $\text{GS}(4)$ provides the second best result. Note, however, that unlike $\text{GS}(4)$ (and ILU preconditioner), our approach does not require one to solve linear problems with large sparse triangular matrices, so it is much faster when parallel architectures are available.

\begin{figure}[!h]
  \centering
    \includegraphics[scale=0.55]{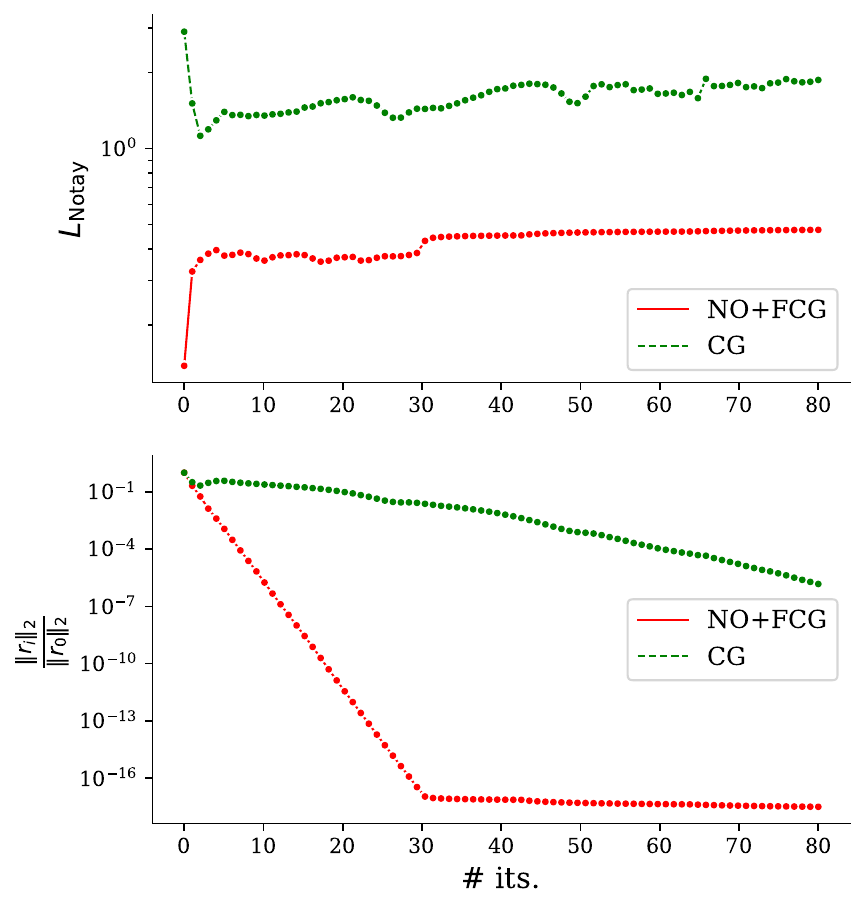}
    \caption{The behavior of $L_{\text{Notay}}$ and the decline of residuals by iteration for Poisson equation with $\text{grid}=32$ in cases of CG and NO+FCG.}
  \label{fig:CG_FCG}
\end{figure}

In addition, we explore how the Notay loss \eqref{eq:final_loss} differs by iteration for two methods: CG and NO+FCG, see \cref{fig:CG_FCG}. We can calculate the $L_{\text{Notay}}$ for CG as it is FCG with $\mathcal{B}=I$. We can see that, for the proposed method, we have a fast convergence and $L_{\text{Notay}} < 1$. It is consistent with the \cref{thm:theorem_1}.  For the CG case, there are  $L_{\text{Notay}} > 1$  and a slow convergence. 

\subsubsection{$L_{\text{Notay}}$ loss vs. $L_2$ loss}
\label{sub:notay_vs_l2}
\begin{table*}[!h]
    \centering
    \begin{tabular}{@{}lccccccc@{}}
        \toprule
        & & \multicolumn{3}{c}{$L_{\text{Notay}}$} & \multicolumn{3}{c}{$L_2$}    \\\cmidrule(r){3-8} 
         & & \multicolumn{3}{c}{$\Vert r_i \Vert_2 / \Vert r_0 \Vert_2$} & \multicolumn{3}{c}{$\Vert r_i \Vert_2 / \Vert r_0 \Vert_2$}  \\\cmidrule(r){3-5} \cmidrule(r){6-8} 
        Dataset & grid &  $10^{-3}$ & $10^{-6}$ & $10^{-12}$& $10^{-3}$ & $10^{-6}$ & $10^{-12}$   \\
        \midrule
        \multirow{3}{*}{Poisson} &  $32$ & 4 & 9 & 20 & 5 & 15 & 34   \\\cmidrule(r){2-5}\cmidrule(r){6-8}
         &  $64$  &  5 & 14 & 31 & 7 & 22 & 53  \\\cmidrule(r){2-5}\cmidrule(r){6-8}
         &  $128$  & 6 & 20 & 48 & 21 & 86 & 210   \\
         \midrule
        \multirow{3}{*}{Diffusion} &  $32$ &  4 & 9 & 31 &  8 & 22 & 50   \\\cmidrule(r){2-5}\cmidrule(r){6-8}
         &  $64$ & 5 & 14 & 36 & 11 & 36 & 80   \\\cmidrule(r){2-5}\cmidrule(r){6-8}
         &  $128$  &5 & 19 & 47 & --- & --- & ---  \\
        
        \bottomrule 
    \end{tabular}
    \caption{The number of iterations for FCG with $\mathcal{B}=\text{NO}$ needed to drop initial residual by three different factors. In table, there are two cases: \textbf{(a)} $\text{NO}$ trained with Notay loss in form \eqref{eq:notay_loss_with_error};
    \textbf{(b)} $\text{NO}$ trained with and with $L_2$-loss in form \eqref{eq:l2_loss_with_error}. In both cases, NO was trained on residuals from Krylov subspace, $r \sim p_{\mathcal{K}_{m}}(r)$.}
    \label{table:notay_l2}
\end{table*}

For this experiment, we have generated training datasets $\mathcal{D}_{\text{train}}$ for Poisson \eqref{eq:poisson} and Diffusion \eqref{eq:diffusion} equations with  $N_{\text{train}} = \text{grid}$. We chose this value because the higher the resolution, the more samples are needed to train NO. To make a comparison of losses, we generated residuals from the Krylov subspace by utilizing $100$ iterations of CG. In total, NO was trained on $100 \cdot N_{\text{train}}$ residuals and errors. The errors were calculated as $e_i^j = A^{-1}_j f_j - u_i^j$, where $ i = 1, \ldots, 100$, ~$j = 1, \ldots, N_{\text{train}}$ and $u_i^j$ is the solution on $i$-th iteration of CG for $j$-th sample. A detailed description of NO and training details are available in \cref{appendix:Architectures and training details}. For NO training, we used Notay and $L_2$ losses in the forms \eqref{eq:notay_loss_with_error} and \eqref{eq:l2_loss_with_error} respectively.

The test datasets $\mathcal{D}_{\text{test}}$ consists of $N_\text{test} = 20$ samples. Then we run $N_{\text{iter}}$ of FCG with trained NO as the preconditioner. We chose $N_{\text{iter}} = 2 \cdot \text{grid}$ because higher resolution requires more iterations to converge. The parameter $m_{\text{max}}$ in \cref{alg:algorithm_1} was selected to be $20$. The numbers of iterations needed to drop the initial residual by factors $10^{3}, 10^{6}, 10^{12}$ are illustrated in \cref{table:notay_l2}.

As can be seen from the \cref{table:notay_l2}, the number of iterations required to reduce residuals by any factor increased when we use $L_2$ loss instead of $L_\text{Notay}$ loss. For factors $10^3$, $10^6$, $10^{12}$, this increase was more than $25\%$, $50\%$, $60\%$, respectively. For $L_\text{Notay}$ loss, the number of iterations for factor $10^{12}$ increases less than $1.55$ times, with a double increase in resolution. For $L_2$ loss, this number of iterations grows more than $1.55$ times with double increase in resolution. Additionally, the missing results for $\text{grid}=128$ means that there was no convergence for this case (see \cref{fig:bad_result_l2} in \cref{appendix:missing values}).  

\begin{figure}[!h]
  \centering
    \includegraphics[scale=0.6]{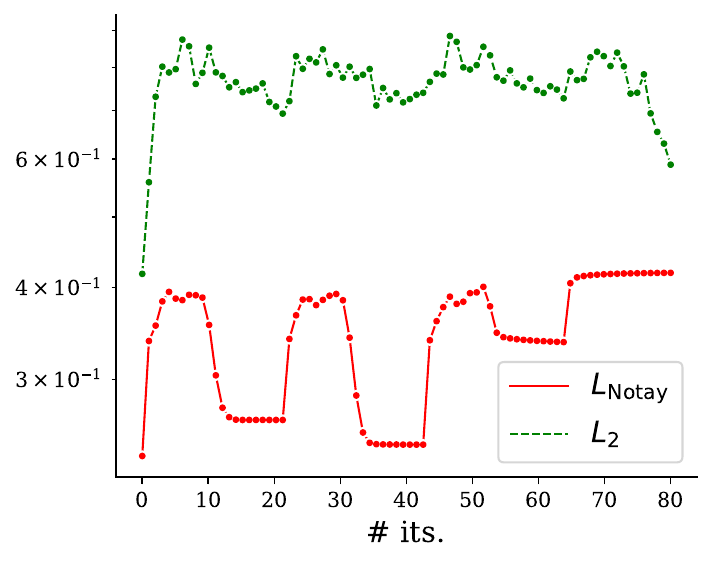}
    \caption{The behavior of $L_{\text{Notay}}$ by iteration for Diffusion equation with $\text{grid}=32$ in cases of NO trained with two different losses ($L_{\text{Notay}}$ and $L_2$).}
  \label{fig:l2_l_notay}
\end{figure}

Furthermore, we explore how the values of \eqref{eq:final_loss} differ by iteration for NO trained with different losses ($L_{\text{Notay}}$ and $L_2$), see \cref{fig:l2_l_notay}.  $L_{\text{Notay}}$ correlates with rate of convergence according to \cref{thm:theorem_1}. We can see that, for both cases, we have $L_{\text{Notay}} < 1$. However, training a NO with the proposed loss \eqref{eq:notay_loss_with_error} gives lower values of function \eqref{eq:final_loss} in all iterations of the FCG than training with a standard $L_2$-loss.

\subsubsection{FCG with NO trained on differently generated residuals}
\begin{table*}[!h]
    \centering
    \begin{tabular}{@{}lccccccc@{}}
        \toprule
        & & \multicolumn{3}{c}{$r \sim p_{\mathcal{K}_{m}}(r)$} & \multicolumn{3}{c}{$r \sim p_{\mathcal{K}_{0}}(r)$}  \\\cmidrule(r){3-8} 
         & & \multicolumn{3}{c}{$\Vert r_i \Vert_2 / \Vert r_0 \Vert_2$} & \multicolumn{3}{c}{$\Vert r_i \Vert_2 / \Vert r_0 \Vert_2$}  \\\cmidrule(r){3-5} \cmidrule(r){6-8}
        Dataset & grid &  $10^{-3}$ & $10^{-6}$ & $10^{-12}$& $10^{-3}$ & $10^{-6}$ & $10^{-12}$   \\
        \midrule
        \multirow{3}{*}{Poisson} &  $32$ & 4 & 9 & 20   & 4 & 10 & 21  \\\cmidrule(r){2-5}\cmidrule(r){6-8}
         &  $64$  &  5 & 14 & 31  &  5 & 18 & 67 \\\cmidrule(r){2-5}\cmidrule(r){6-8}
         &  $128$  & 6 & 20 & 48  & 7 & 49 & 153 \\
         \midrule
        \multirow{3}{*}{Diffusion} &  $32$ &  4 & 9 & 31  &  5 & 23 & 56 \\\cmidrule(r){2-5}\cmidrule(r){6-8}
         &  $64$ & 5 & 14 & 36  & 6 & --- & ---  \\\cmidrule(r){2-5}\cmidrule(r){6-8}
         &  $128$  &5 & 19 & 47 & 10 & --- & ---  \\
        
        \bottomrule 
    \end{tabular}
    \caption{The number of iterations for FCG with $\mathcal{B}=\text{NO}$ needed to drop initial residual by three different factors. In table, there are two cases: \textbf{(a)} $\text{NO}$ trained on residuals from Krylov subspace;
    \textbf{(b)} $\text{NO}$ trained on residuals obtained from random right-hand sides. In both cases, NO was trained with Notay loss in the form \eqref{eq:notay_loss_with_error}.}
    \label{table:notay_krylov}
\end{table*}
Train datasets $\mathcal{D}_{\text{train}}$ for both equations had $N_{\text{train}} = \text{grid}$. For residuals $r \sim p_{\mathcal{K}_{m}}(r)$, the training procedure is the same as in \cref{sub:notay_vs_l2}. For random residuals, we generated $100 \cdot N_{\text{train}}$ samples of $u_i^j \sim \mathcal{N}(0, I)$ and $f_i^j \sim p_f$, where $i = 1, \ldots, 100$, ~$j=1, \ldots, N_{\text{train}}$. Therefore, the residuals and errors were calculated as follows: $r_i^j = f_i^j - A_j u_i^j$,  $e_i^j = A^{-1}_j f_i^j - u_i^j$. We trained NO with Notay loss using the architecture and training setup described in \cref{appendix:Architectures and training details}.

The test datasets $\mathcal{D}_{\text{test}}$ contain $N_{\text{test}} = 20$ samples. The FCG algorithm has the same number of iterations and the same parameter $m_{\text{max}}$ as in \cref{sub:notay_vs_l2}. The results are shown in \cref{table:notay_krylov}.

It can be seen from the data in \cref{table:notay_krylov} that for the simplest case (Poisson equation, $\text{grid}=32$), we got almost the same results for both types of residuals. However, with increasing resolution, there was a significant deterioration in convergence for residuals not from the Krylov subspace. For the Poisson equation, the number of iterations for factor $10^{12}$ increased more than $2$ times when the resolution doubled. For the Diffusion equation, the algorithm stopped converging when the grid is larger than $32$ (see \cref{fig:bad_result_random_64}, \cref{fig:bad_result_random_128} in \cref{appendix:missing values}).

\subsection{Different grid approach}
\begin{table*}[!h]
    \centering
    \begin{tabular}{@{}lccccc@{}}
        \toprule
         & &  & \multicolumn{3}{c}{$\Vert r_i \Vert_2 / \Vert r_0 \Vert_2$} \\\cmidrule(r){4-6}
        Dataset & Train grid & Test grid &  $10^{-3}$ & $10^{-6}$ & $10^{-12}$ \\
        \midrule
        \multirow{4}{*}{Poisson} &  \multirow{1}{*}{$32$} &  \multirow{1}{*}{$64$}  & 9 & 18 & 39  \\\cmidrule(r){2-6}
         &  \multirow{1}{*}{$64$} &  \multirow{1}{*}{$128$}  & 13 & 31 & 72 \\\cmidrule(r){2-6}
         &  \multirow{1}{*}{$32$} &  \multirow{1}{*}{$128$}  & 30 & 79 & 191  \\\cmidrule(r){2-6}
         &  \multirow{1}{*}{$64$} &  \multirow{1}{*}{$256$}  & 40 & 104 & 245  \\
         \midrule
        \multirow{4}{*}{Diffusion} &  \multirow{1}{*}{$32$} &  \multirow{1}{*}{$64$}  & 8 & 20 & 52 \\\cmidrule(r){2-6}
         &  \multirow{1}{*}{$64$} &  \multirow{1}{*}{$128$}  & 8 & 22 & 50  \\\cmidrule(r){2-6}
         &  \multirow{1}{*}{$32$} &  \multirow{1}{*}{$128$}  &  17 & 59 & 170 \\\cmidrule(r){2-6}
         &  \multirow{1}{*}{$64$} &  \multirow{1}{*}{$256$}  & 27 & 133 & ---   \\
        
        \bottomrule 
    \end{tabular}
    \caption{Number of iterations for FCG with $\mathcal{B}=\text{NO}$ needed to drop initial residual by three different factors. The table consists of results for different grids on training and testing.}
    \label{table:results_diff_grids}
\end{table*}
The testing process becomes more complex when conducted on a higher resolution than the one used during training. Therefore, we trained NO on more samples than in previous experiments. There were $N_{\text{train}} = 100$ samples in the training datasets $\mathcal{D}_{\text{train}}$. For generation residuals and errors, we used $100$ iterations in CG. The details of the architecture and the training process are presented in \cref{appendix:Architectures and training details}.

We tested FCG with NO on $N_{\text{test}} = 20$ samples. The resolutions of the test datasets $\mathcal{D}_{\text{test}}$ were twice or four times higher than the resolutions of the train datasets $\mathcal{D}_{\text{train}}$. Thus, FCG needed more iterations to converge. We selected $N_{\text{iter}} = 300$ for all cases. The results are demonstrated in \cref{table:results_diff_grids}.

The results from \cref{table:results_diff_grids} demonstrate that the operator can serve as an effective preconditioner even if it is trained at a lower resolution. For the case of missing results, there were not enough iterations set for the FCG algorithm (see \cref{fig:bad_result_diff_grid} in \cref{appendix:missing values}). 
\label{section:experiments}

\section{Related research}
\label{section:related research}
Literature on constructing preconditioners is vast, so we do not attempt a systematic review. In place of that, we will mention several contributions that allow us to frame our research properly.

As a first class of work, we would like to mention methods that construct preconditioners by solving auxiliary optimization problems (e.g., loss \eqref{eq:final_loss} that we used). Typically, such works restrict somehow the form of matrix $M$ used as a preconditioner and optimize $\left\|I - M^{-1}A\right\|$ or $\left\|M - A\right\|$. For example, in \cite{tyrtyshnikov1992optimal} author showed how to efficiently construct general circulant preconditioners, and in \cite{grote1997parallel} authors exploited sparsity constraints.

A similar strategy is widely used with models based on neural networks. For example, in \cite{cui2022fourier} authors train a neural network to output eigenvalues of preconditioner in Fourier basis in effect reproducing the circulant preconditioner. Similarly, authors of \cite{hausner2023neural} utilize a graph neural network with sparsity constraints on Cholesky factors to train linear preconditioner for CG.

The second class of works is learning generalized preconditioners, i.e., matrix $M$ used to transform the linear equation $MAx=Mb$ to improve the convergence of base iterations (Richardson, Jacobi, multigrid). As an example of such approaches, we can mention \cite{hsieh2019learning} and \cite{zhang2022hybrid}. In the first contribution authors train the linear U-Net model as a preconditioner for Richardson iterations, in the second contribution authors pursue a similar idea with DeepONet architecture but with Jacobi and multigrid.

Compared with the mentioned works, our approach is more direct since FCG, which we use as a base iterative method, is not restricted to working with linear preconditioners. Besides, we directly learn nonlinear operators, and with that, we avoid the need to consider structured matrices that guarantee cheap matrix-vector products. Finally, the proposed loss function is markedly different from what is used in other contributions.

\section{Conclusion}
In the present research, we have suggested using a trained neural operator as a nonlinear preconditioner for the flexible conjugate gradient method. The results of our study demonstrated the superiority of the proposed approach over utilizing different classical preconditioners. In addition, we introduced a novel loss function derived from the energy norm that guarantees convergence and achieves superior results compared to the $L_2$ loss function used in training. Moreover, we have implemented a novel learning scheme incorporating random vectors derived from the Krylov subspace. Based on the results of the ablation study, the utilization of classical random right-hand sides instead of the Krylov subspace severely diminished the effectiveness of training. Our findings also demonstrated that neural operators can be trained at lower resolutions and effectively act as a cost-effective preconditioner for higher resolutions. 

The limitation of the proposed method is that the matrix $A$ should be symmetric positive definite as we use CG. Our potential for utilizing our approach on higher resolutions is limited by GPU memory. 

In the future, it is planned to explore the proposed approach with other types of neural operators. Additionally, in this paper, we only focused on using a uniform mesh. By utilizing different polynomials in SNO, we can extend our approach to non-uniform grids.  Moreover, we plan to explore other groups of PDEs for which the large sparse linear algebraic system has a symmetric positive definite matrix.

\section{Impact Statement}
This paper presents work whose goal is to advance the field of Machine Learning. There are many potential societal consequences of our work, none which we feel must be specifically highlighted here.


\bibliography{FCG/refs}
\bibliographystyle{icml2024}

\newpage
\appendix
\onecolumn
\section{Architectures and training details.}
\label{appendix:Architectures and training details}
In this section, we provide details of architecture of used NO and training process. 

We use SNO in Fourier basis (see \cite{fanaskov2022spectral}) with encoder-processor-decoder architecture. Number of features in the processor for Poisson equation is $32$. For Diffusion equation, we change number of features according to grid, see \cref{table:n_params}. Number of SNO layers is $4$ and number of orthogonal polynomials is $20$. We utilize GeLU as an activation function. 
\begin{table*}[!ht]
    \centering
    \begin{tabular}{@{}lcccc@{}}
        \toprule
        \multirow{2}{*}{Diffusion} & grid & $32$ & $64$ & $128$ \\\cmidrule(r){2-5}
         & $N_{\text{params}}$ & $32$ & $48$ & $85$ \\
        \bottomrule 
    \end{tabular}
    \caption{Number of parameters in processor of SNO.}
    \label{table:n_params}
\end{table*}

\begin{table}[h]
    \centering
    \vskip 0.15in
    \begin{tabular}{@{}lcccccc}
        \toprule
        Dataset & grid & $\nu$ & $\nu$ decay / epoch & weight decay & $N_{\text{epoch}}$ & $N_{\text{batch}}$  \\ \midrule
        \multirow{2}{*}{Poisson} & $32$, $64$ & $10^{-3}$ & $0.5 \big/ 50$  &  $10^{-2}$ &  $200$ & $32$ \\\cmidrule(r){2-7}
         & $128$ &  $10^{-3}$ &  $0.5 \big/ 50$ & $10^{-2}$ & $200$ &  $8$\\ \midrule
        \multirow{2}{*}{Diffusion} & $32$, $64$ & $10^{-3}$  & $0.5 \big/ 50$ & $10^{-2}$ & $150$ & $16$  \\\cmidrule(r){2-7}
         & $128$ &  $10^{-3}$ & $0.5 \big/ 50$ & $10^{-2}$ & $150$ & $4$\\
    \bottomrule
    \end{tabular}
        \caption{Training details: $\nu$ --- learning rate, $\nu$ decay / epoch --- weight decay per epoch, $N_{\text{epoch}}$ --- number of epoch used for training, $N_{\text{batch}}$ --- batch size.}
    \label{table:training_details}
    \vskip -0.1in
\end{table}

\section{FCG convergence for experiments with missing values in results}
\label{appendix:missing values}
\begin{figure}[!h]
  \centering
    \includegraphics[scale=0.66]{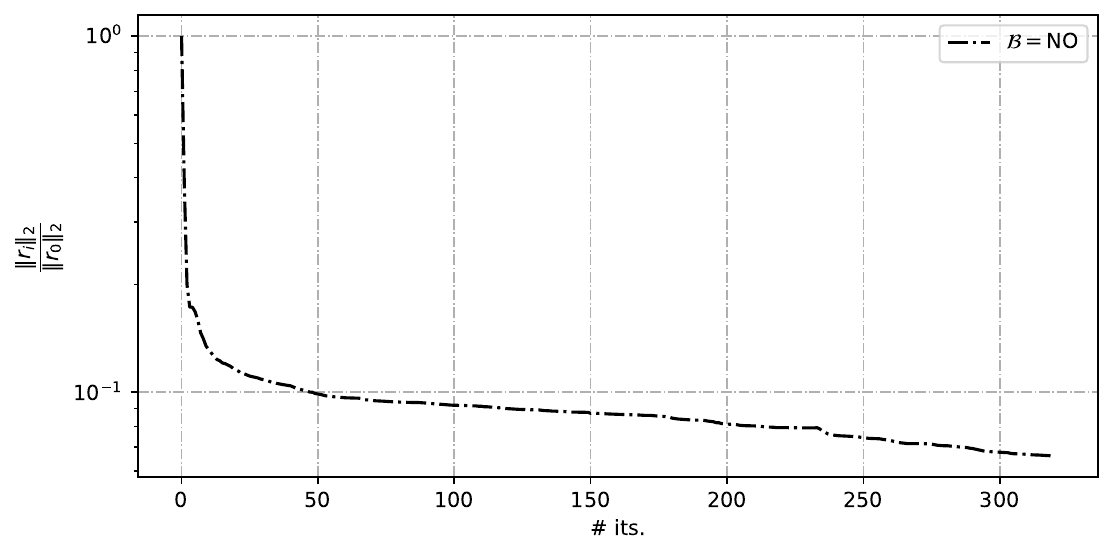}
    \caption{The decline of residuals by iteration for Diffusion equation with $\text{grid}=128$ in $L_2$, $r \sim p_{\mathcal{K}_{m}}(r)$ case.}  
  \label{fig:bad_result_l2}
\end{figure}

\begin{figure}[!h]
  \centering
    \includegraphics[scale=0.66]{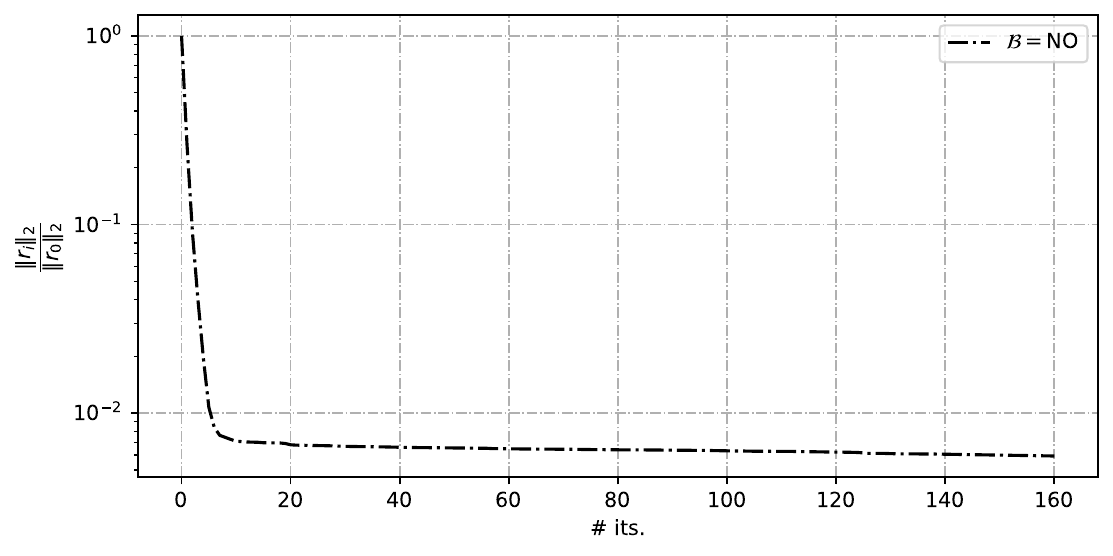}
    \caption{The decline of residuals by iteration for Diffusion equation with $\text{grid}=64$ in $L_{\text{Notay}}$, $r \sim p_{\mathcal{K}_{0}}(r)$ case.}
  \label{fig:bad_result_random_64}
\end{figure}

\begin{figure}[!h]
  \centering
    \includegraphics[scale=0.66]{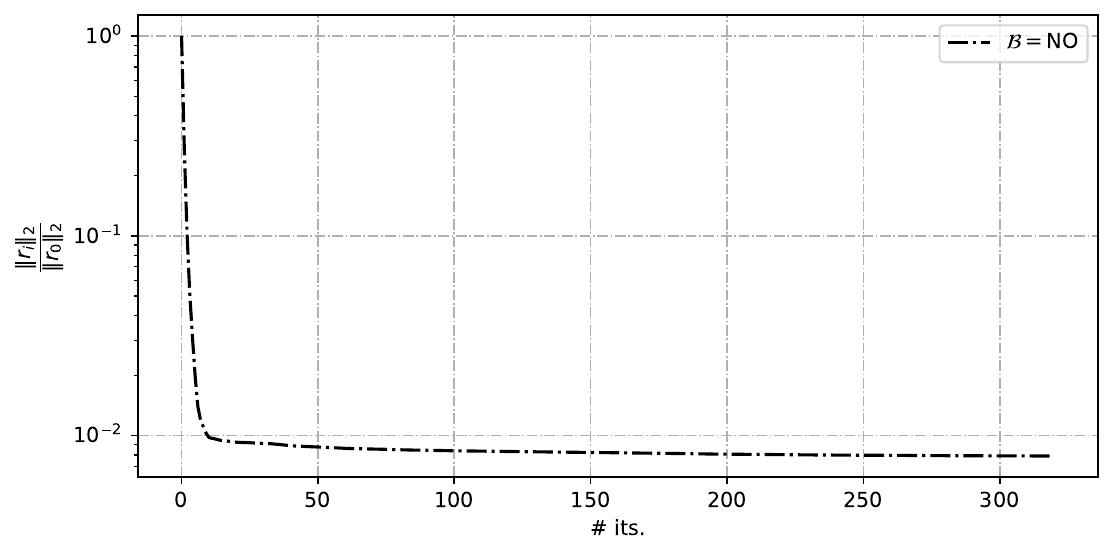}
    \caption{The decline of residuals by iteration for Diffusion equation with $\text{grid}=128$ in $L_{\text{Notay}}$, $r \sim p_{\mathcal{K}_{0}}(r)$ case.}
  \label{fig:bad_result_random_128}
\end{figure}

\begin{figure}[!h]
  \centering
    \includegraphics[scale=0.66]{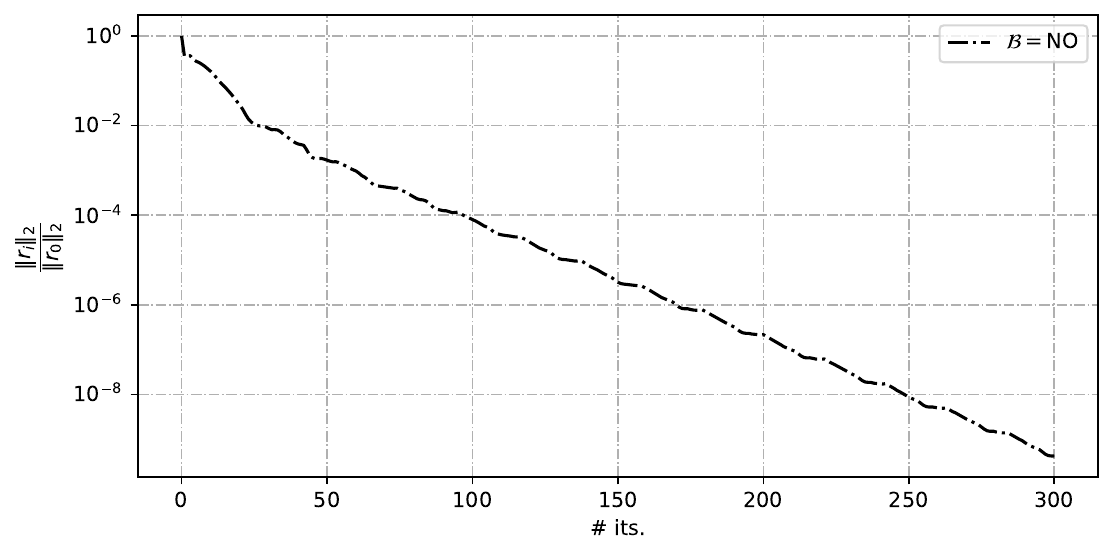}
    \caption{The decline of residuals by iteration for  Diffusion equation with $\text{grid}_\text{train}=64$ and $\text{grid}_\text{test}=256$ in different grids case.}
  \label{fig:bad_result_diff_grid}
\end{figure}

\end{document}